\documentstyle[12pt,amsfonts,epsfig]{article}

\newcommand{\new}{\newcommand}

\newcommand\hdiff{{\cal H}_{\rm diff}}
\new{\qed}{{\hfill \rule{1ex}{1ex}}\vskip 1em}
\new{\et}{\hspace{-0.08in}{\bf .}\hspace{0.1in}}
\new{\tnsr}{\otimes} \new{\iso}{\cong} \new{\union}{\cup}
\new{\implies}{\Rightarrow} \new{\maps}{\colon} \new{\A}{{\cal A}}
\new{\G}{{\cal G}} \new{\R}{{\Bbb R}} \new{\C}{{\Bbb C}}
\new{\tr}{{\rm tr}} \new{\range}{{\rm Range}} \new{\dom}{{\rm Domain}}
\new{\Lie}{{\rm Lie}} \new{\Fun}{{\rm Fun}} \new{\Aut}{{\rm Aut}}
\new{\sgn}{{\rm sgn}}
\new{\SU}{{\rm SU}}
\new{\U}{{\rm U}} \new{\abs}[1]{\left|#1 \right|}
\new{\norm}[1]{\left\| #1 \right\|} \new{\bracket}[1]{\langle #1
\rangle} \new{\defequals}{\stackrel{\rm def}{=}}
\new{\into}{\hookrightarrow} \new{\comb}[2]{{#1 \choose #2}}
\new{\lbl}[1]{\label{#1} \if\draft y
\smash{\makebox[0pt]{\hspace{-0.5in} \raisebox{8pt}{\rm\tiny #1}}} \fi}

\new{\pic}[5]{\raisebox{#3pt}{
\hspace{#4pt}\epsfig{file=#1.eps,height=#2pt}\hspace{#5pt}}}

\newtheorem{proposition}{Proposition}
\newtheorem{theorem}{Theorem}

\newcounter{letter} \newcounter{numeral} \newcounter{Numeral}
\newenvironment{alphalist}{
\begin{list}{{(\alph{letter})}}{\usecounter{letter}}
}{\end{list}}

\newenvironment{pf}{{\noindent\bf Proof:} }{\qed}

\begin{document}

      \begin{center}
      {\bf Diffeomorphism-Invariant Spin Network States \\}
      \vspace{0.5cm}
      {\em John C. Baez\\}
      \vspace{0.3cm}
      {\small  Department of Mathematics \\
      University of California\\
      Riverside, CA 92521-0135\\}
      \vspace{0.3cm}
      \renewcommand{\thefootnote}{\fnsymbol{footnote}}
      {\em Stephen Sawin \\}
      \vspace{0.3cm}
      {\small  Department of Math and C.S. \\
      Fairfield University \\
        Fairfield, 06430-5195 \\ }
      \vspace{0.3 cm}
      {\small email: \tt baez@math.ucr.edu, ssawin@fair1.fairfield.edu}
      \end{center}

\begin{abstract} 
We extend the theory of diffeomorphism-invariant spin network states
from the real-analytic category to the smooth category.  Suppose that 
$G$ is a compact connected semisimple Lie group and $P \to M$ is a smooth
principal $G$-bundle.  A `cylinder function' on the space of
smooth connections on $P$ is a continuous complex function of the
holonomies along finitely many piecewise smoothly immersed curves in
$M$.  We construct diffeomorphism-invariant functionals on the space
of cylinder functions from `spin networks': graphs in $M$ with edges 
labeled by representations of $G$ and vertices labeled by intertwining 
operators.  Using the `group
averaging' technique of Ashtekar, Marolf, Mour\~ao and Thiemann, we
equip the space spanned by these `diffeomorphism-invariant spin
network states' with a natural inner product.
\end{abstract}

\section*{Introduction}

In the `new variables' approach to quantizing gravity, the kinematical
Hibert space of the theory should consist of functions on some
completion of the space of connections on a principal $\SU(2)$ bundle
over the smooth 3-manifold representing space.  Defining the inner
product in this Hilbert space requires a measure on the completed space
of connections.  Starting from this Hilbert space, one should then solve
constraint equations corresponding to gauge-invariance,
diffeomorphism-invariance and invariance under time evolution to obtain
the space of physical states.  The deepest problems with this program
are those associated with invariance under time evolution, i.e.\ the
Hamiltonian constraint.  However, the proper treatment of the other two
constraints also presents problems, some of which have been avoided by
assuming $M$ is real-analytic and only demanding invariance under
real-analytic diffeomorphisms.  Working in the real-analytic category,
Ashtekar, Lewandowski, Marolf, Mour\~ao and Thiemann have constructed a
Hilbert space of gauge- and diffeomorphism-invariant states spanned by
`spin networks' \cite{ALMMT}.  Here we do the same in the smooth
category.

Before describing our result more precisely, let us briefly review the
state of the art \cite{AI,AL,B,BS,RS}.  First, there is a
natural way to complete the space of connections on any smooth
principal $G$-bundle $P \to M$ when $G$ is a compact connected Lie
group.  This goes as follows.  Let $\A$ be the space of smooth
connections on $P$, and define a `cylinder function' on $\A$ to be a
continuous function of the holonomies along finitely many piecewise
immersed paths in $M$.  Taking the $\sup$ norm completion of the
algebra of cylinder functions, one obtains a commutative C*-algebra,
and the spectrum of this C*-algebra is a compact Hausdorff space
$\overline\A$ having $\A$ as a dense subset.   Any cylinder function
uniquely extends to a continuous function on $\overline \A$.

Second, there is a natural Borel measure on $\overline\A$, the `uniform
measure' $\mu_0$.  With respect to this measure, the probability
distribution of the holonomies along the edges of a graph embedded in
$M$ is given by Haar measure on a product of copies of $G$.  In other
words, suppose that $e_1,\dots,e_n$ are smoothly embedded  copies of the
unit interval in $M$ which intersect, if at all, only at their
endpoints.  In this situation we say that $e_i$ are the edges of a
graph.  Let $F$ be the cylinder function
\[   F(A) = f({\rm T}\exp \int_{e_1}A,\dots,{\rm T}\exp \int_{e_n}A)\]
where $f$ is a continuous complex-valued function on $G^n$ and ${\rm T}
\exp\int_{e_i}A$ is the holonomy of the connection $A$ along $e_i$,
regarded as a group element by means of an arbitrary trivialization of
$P$ at the endpoints of this curve.  Then the integral of $F$ with
respect to the uniform measure is given by
\[     \int_{\overline \A} F(A)\, d\mu_0(A)
= \int_{G^n} f(g_1,\dots,g_n) \,dg_1 \cdots dg_n \]
where the right-hand integral is taken with respect to normalized
Haar measure on $G^n$.

In the real-analytic category the above property is sufficient to
characterize the uniform measure, because any cylindrical function can
be expressed as above in terms of the holonomies of curves forming a 
graph.  This no longer holds in the smooth category, making it a bit
trickier to fully characterize uniform measure.  The reason is that
smoothly embedded curves can intersect each other in extremely
complicated ways, for example in a Cantor set.   In a previous paper
\cite{BS} we dealt with this issue by generalizing graphs to `webs'.
Like a graph, a web consists of finitely many curves embedded in $M$,
but in a web these curves are allowed to intersect each other in certain
specified ways. In Section 1 we recall the concept of a web and
characterize uniform measure in the smooth category using webs.

Using uniform measure one can define the `kinematical Hilbert space'
$L^2(\overline\A)$.  The group of gauge transformations has a unitary
representation on this Hilbert space, coming from its action on $\A$,
and there is a large subspace of $L^2(\overline \A)$ consisting of 
functions invariant under all gauge transformations.   We denote this
`gauge-invariant Hilbert space' by $L^2(\overline{\A/\G})$, since its
elements may also be regarded as square-integrable functions on a
certain completion of the space of connections modulo gauge
transformations.   In our previous paper we constructed an explicit set
of functions spanning the gauge-invariant Hilbert space, the `spin web
states'. Each such state is determined by a `spin web': a web with edges
labeled by representations of $G$ and vertices labeled by intertwining
operators.  We review the theory of spin webs in Section 2.

The simplest spin web states are the `spin network states',
corresponding to webs that are simply graphs.  In the real-analytic
category, spin network states are enough to span the gauge-invariant
Hilbert space.  This is no longer true in the smooth category.
However, we show in Section 2 that any spin network state is
orthogonal to any spin web state that cannot be written as a spin
network state with the same underlying graph.  This result is crucial
for the next step: constructing the diffeomorphism-invariant Hilbert
space.

Naively one might try to define the diffeomorphism-invariant Hilbert
space as the subspace of $L^2(\overline{\A/\G})$ consisting of
functions invariant under all diffeomorphisms of $M$, or at least
those in the identity component of the diffeomorphism group.
Moreover, one might naively be inclined to obtain such functions by
`averaging over the action of the diffeomorphism group'.  However,
things are not so simple: there appears to be no `Haar measure' on the
diffeomorphism group, and there are typically very few
diffeomorphism-invariant functions in $L^2(\overline{\A/\G})$.  The
point is that one should seek diffeomorphism-invariant elements, not
of $L^2(\overline{\A/\G})$, but of some larger space of `generalized
functions' on $\overline{\A/\G}$.  Working in the real-analytic
category, Ashtekar, Lewandowski, Marolf, Mour\~ao and Thiemann
\cite{ALMMT} constructed such diffeomorphism-invariant generalized
functions by a clever procedure which amounts to averaging 
spin network states over the action of the diffeomorphism group.

Using the orthogonality result of the previous section, in Section 3 we
carry out a similar group averaging procedure in the smooth category,
obtaining `diffeomorphism-invariant spin network states' labeled by
diffeomorphism equivalence classes of spin networks.  Completing the
space spanned by these in its natural inner product, we obtain the
diffeomorphism-invariant Hilbert space.

One might imagine extending this diffeomorphism-invariant Hilbert space
to include more general diffeomorphism-invariant spin web states.
Unfortunately, spin webs that are not spin networks behave badly.  In
general two spin webs may not be orthogonal even if their underlying
webs do not have the same range.  Further, there are infinitely many
diffeomorphisms taking a typical spin web state to spin web states that
are not orthogonal to it, even after the obvious quotients.  We give
examples of both these phenomena in Section 4.  Thus it appears difficult
to find an explicit orthonormal basis of the full gauge-invariant
Hilbert space, and difficult to construct diffeomorphism-invariant
states from spin webs.  

We should note that when $M$ is real-analytic, the previously studied
real-analytic spin networks are a special case of our smooth spin
networks.  Furthermore, any smooth manifold can be given a
real-analytic structure, which is unique up to smooth diffeomorphism.
However, there are many more smooth spin networks than real-analytic
ones, even modulo smooth diffeomorphisms, because there are smooth
vertex types unrealizable by analytic curves.

\section{Uniform measure}

We begin with a terse review of uniform measure on the space of
connections.   For the most part we follow the treatment in our previous
paper \cite{BS}, but we simplify the setup using a result of Lewandowski
and Thiemann \cite{LT}.

Fix a connected compact semisimple Lie group $G$, a smooth (paracompact)
manifold $M$, and a smooth principal $G$-bundle $P \to M$.   By a {\em
curve} we mean a piecewise smooth map from an interval $[a,b] \subset
\R$ to $M$ that is an immersion on each piece.  Two curves are
considered equivalent if one is the composition of the other with an
orientation-preserving diffeomorphism between their domains (so that
one is just a reparametrization on the other).  A {\em family} is a
finite set of curves with a chosen ordering $c_1, \ldots ,c_n$.  If
$C$ is such a family, we define $\range(C)$, the {\em range\/} of $C$,
to be the union of the ranges of the individual curves.  A point $p$
in $\range(C)$ is a {\em regular point\/} if it is not the image of an
endpoint or nondifferentiable point of $C$, and there is a
neighborhood of it whose intersection with $\range(C)$ is an embedded
interval.  A family $C$ is {\em parametrized consistently\/} if each
curve is parametrized so that $c_i(t)=c_j(s)$ implies $t=s$.  Thus
each of the curves is actually an embedding, and each point $p$ in the
range of the family is associated to a unique value of the parameter,
which we call $t(p)$.  If a family $\{c_1, \ldots, c_n\}$ is
parametrized consistently and $p$ is a point in $\range(C)$, define
the {\em type\/} of a regular point $p$, $\tau_p$, to be the Lie
subgroup of $G^n$ consisting of all $n$-tuples $(g_1, \dots, g_n)$
such that for some $g \in G$ we have $g_i = g$ if $p$ lies on $c_i$,
and $g_i = 1$ otherwise.

A family $T$ is a {\em tassel based on\/} $p \in \range(T)$
if:
\begin{alphalist}
\item $\range(T)$ lies in a contractible open subset of $M$.
\item $T$ can be consistently parametrized in such a way that $c_i(0)
= p$ is the left endpoint of every curve $c_i$.  
\item Two curves in $T$ that intersect at a point other than $p$.
intersect at a point other than $p$ in every neighborhood of $p$.
\item Any type which occurs at some point in $\range(T)$ occurs in
every neighborhood of $p$.
\item No two curves in $T$ have the same range.
\end{alphalist}
A {\em web\/} $\omega$ is a finite collection of tassels 
$\omega^1, \dots, \omega^k$ such that for $i \ne j$:
\begin{alphalist}
\item Any curve in the tassel $\omega^i$
intersects any curve in $\omega^j$, if at all, only at their endpoints
\item There is a neighborhood of each such
intersection point whose intersection with $\range(\omega^i \cup
\omega^j)$ is an embedded interval.
\item $\range(\omega^i)$ does not contain the base of $\omega^j$.
\end{alphalist}
We may apply concepts defined for families to webs, since every web
$\omega$ has an associated family $\omega^1 \cup \cdots \cup \omega^k$.
We define an {\em edge\/} of a web $\omega$ to be a curve in one of the
families $\omega^i$, and define a {\em vertex\/} of $\omega$ to be a
point of $M$ that is an endpoint of some edge of $\omega$.  Since the
edges of a web are oriented, we may speak of the {\em source\/} and {\em
target\/} of any edge, these being its initial and final endpoints.

Using webs one can characterize uniform measure as follows.
First, note that any web $\omega$ with edges
$e_1,\dots,e_n$ together with trivializations of $P$ at the vertices
of $\omega$ determines a map from $\A$ to $G^n$ given by:
\[   A \mapsto ({\rm T}\exp \int_{e_1}A,\dots,{\rm T}\exp
\int_{e_n}A).\]
This map extends uniquely to a continuous map
\[   p_\omega  \maps \overline \A \to G^n .\]
We may push forward any Borel measure on $\overline \A$
to a Borel measure on $G^n$ by this map $p_\omega$.   We then have:

\begin{proposition}\et  There exists a unique Borel measure on
$\overline
\A$, the {\rm uniform measure} $\mu_0$, such that if $\omega$ is a web
with $n$ edges, the pushfoward of $\mu_0$ by $p_\omega$ is normalized
Haar measure on $G^n$.  \end{proposition}

\begin{pf} In our previous paper we showed that the range of $p_\omega$
is a Lie subgroup of $G^n$ and that $\mu_0$ is uniquely characterized by
the property that its pushforward is normalized Haar measure on this
subgroup for every web $\omega$.  Previously we did not include clause
(e) in the definition of a tassel, but we may assume this without loss
of generality, since two consistently oriented curves with the same
range have the same holonomy for every connection.   Using this clause
and the fact that $G$ is semisimple, Lewandowski and Thiemann
\cite{LT} subsequently showed that $p_\omega$ is onto.  \end{pf}

\section{Spin webs and spin networks}

Using uniform measure one can define the kinematical Hilbert space
$L^2(\overline\A)$.  This in turn allows us to construct the
gauge-invariant Hilbert space $L^2(\overline{\A/\G})$, which consists
of all functions in $L^2(\overline\A)$ that are invariant under gauge
transformations.  In this section we provide a detailed description of
the gauge-invariant Hilbert space in terms of spin webs and spin
networks.

A {\em spin web\/} is a triple  $W = (\omega,\rho,\iota)$ consisting of:
\begin{alphalist}
\item a web $\omega$
\item a labeling $\rho$ of each edge $e$ of $\omega$ with a
nontrivial irreducible representation $\rho_e$ of $G$
\item a labeling
$\iota$ of each vertex $v$ of $\omega$ with an intertwining operator
$\iota_v$ from the tensor product of the $\rho_e$ for which $v$ is
the target of $e$ to the tensor product of the $\rho_e$ for which
$v$ is the source of $e$.
\end{alphalist}
Given a spin web $W = (\omega,\rho,\iota)$, the {\em spin web state\/}
$\Psi_W$ is the cylinder function on $\A$ constructed as follows:
\[
\Psi_W(A) = \bigl [\bigotimes_e\,
\rho_e({\rm T}\exp \int_{e} A)\bigr ]\cdot\bigl
[\bigotimes_v\,\iota_v \bigr ],
\]
where `$\cdot$' stands for contracting, at each vertex $v$ of $\omega$,
the upper indices of the matrices corresponding to the incoming edges,
the lower indices of the matrices assigned to the outgoing edges, and
the corresponding indices of the intertwiner $\iota_v$.

\begin{proposition}\et Finite linear combinations of spin web states
are dense in $L^2(\overline{\A/\G})$.  \end{proposition}

\begin{pf}  This is a slight rephrasing of a result in our previous
paper.   There we called every state of the form $\Psi_W$ a `spin
network  state', but here we reserve that term for a special case (see
below).  Also, here we assume without loss of generality that the
underlying web of $W$ has no two edges with the same range and no edge
labeled by a trivial representation.  \end{pf}

We say a subset of $M$ is a {\em graph in $M$\/} if it is the union of a
finite collection of embedded copies of the unit interval that
intersect, if at all, only at their endpoints.  We say that a web is a
{\em graph\/} if its range is a graph in $M$.   Note that for every
graph $G$ in $M$ there is a web with $G$ as its range, and this web is
unique up to inserting and deleting  bivalent vertices and reversing
orientations of edges.

We define a {\em spin network\/} to be a spin web $\Gamma =
(\gamma,\rho,\iota)$ whose underlying web $\gamma$ is a graph.  In this
case we call the spin web state $\Psi_\Gamma$ a {\it spin network
state}. This definition of spin network state is a bit different from
the usual one \cite{B}.  However, apart from the fact that our graphs
have smooth rather than real-analytic edges, the differences are purely
superficial.   To see this, suppose we have a spin network $\Gamma =
(\gamma,\rho,\iota)$ as defined above.  Then the range of $\gamma$ is
the union of the ranges of a finite set of curves intersecting only at
their endpoints.  Subdividing these curves if necessary, we may assume
that each edge of $\omega$ is a product of these curves and their
inverses.   Call this set of curves $E$ and the set of their endpoints
$V$.  Then, just as in the usual definition of spin network state, we
may write
\[
\Psi_\Gamma(A) =
\bigl [\bigotimes_{e\in E}\, \rho_e({\rm T}\exp \int_{e} A)\bigr ]
\cdot\bigl [\bigotimes_{v\in V}\,\iota_v \bigr ],
\]
for some choice of representations $\rho_e$ and intertwining
operators $\iota_v$.

We say that a spin web state $\Psi$ is {\em supported\/} on a web
$\omega$ if it equals $\Psi_W$ for some spin web $W$ with $\omega$ as
its underlying web.  We also say that $\Psi$ is supported on
the range of $\omega$, especially when $W$ is a spin network, so that
the range of $\omega$ is a graph in $M$.   Note that many different
graphs $\gamma$ have the same graph in $M$ as their range: we can 
change $\gamma$ without changing its range by introducing and deleting 
bivalent vertices on embedded intervals, and also by reversing the
orientation of edges.  If a spin network state is supported on $\gamma$ 
it is also supported on all other graphs with the same range.  However,
it is supported on a unique graph in $M$.  The following is a restatement
of arguments in earlier work on spin networks \cite{B}.

\begin{proposition}\et
If two spin network states have nonzero inner product, they are
supported on the same graph in $M$.  
\end{proposition}

\begin{pf}
Let $\Gamma$ and $\Gamma'$ be spin networks with underlying
graphs $\gamma$ and $\gamma'$, respectively.   If the ranges of $\gamma$
and $\gamma'$ are not the same, consider an open segment of an edge
contained in one and not in the other.   We may choose a web on which
both $\Psi_\Gamma$ and $\Psi_{\Gamma'}$ are supported \cite{BS}, and
having one edge lying entirely in this open segment.  Uniform measure
for this web gives the holonomy of this edge Haar measure distribution,
independent as a random variable from the holonomies of the other edges.
In computing the inner product of $\Psi_\Gamma$ and $\Psi_{\Gamma'}$,
this variable will appear once, represented in the nontrivial
irreducible representation labeling that edge.  Since the integral
against Haar measure of the nontrivial irreducible representation of a
group-valued variable is zero, the whole inner product is zero.
Two spin network states with nonzero inner product must therefore
be supported on the same graph $G$ in $M$. 
\end{pf}

The above proposition gives an essentially complete description of the
inner product on the portion of $L^2(\overline{\A/\G})$ spanned by spin
network states, since the description of the inner product of two spin
network states supported on the same graph is well-understood and
involves only elementary group representation theory \cite{B}.   It
remains to understand the inner product of a spin network state with a
general spin web state and the inner product of two arbitrary spin web
states.  The latter question is quite subtle and appears to admit no
simple answer (see Section 4).  The former proves to be tractable, and
is the subject of the next theorem, the key technical result of this
paper.

\begin{theorem} \et
If the inner product of a spin network state with a spin web state is
nonzero, then they are both spin network states supported on the same
graph in $M$.
\end{theorem}

\begin{pf}
Let $\Gamma$ be a spin network with underlying graph $\gamma$, and $W$ a
spin web with underlying web $\omega$.    We assume the inner product of
$\Psi_\Gamma$ and $\Psi_W$ is nonzero and show that the base of any
tassel $\omega^i$ of $\omega$ has a neighborhood $N$ such that
$\range(\omega^i) \cap N \subseteq \range(\gamma)$.   This implies that the
range of $\omega$ is a graph in $M$, so that $\Psi_W$ is a spin network 
state.  The rest of the theorem follows from Proposition 3.

To compute the inner product of $\Psi_\Gamma$ and $\Psi_W$, note from
our previous paper that there is a web $\omega'$ such that every curve
in $\omega$ or $\gamma$ is a product of curves in $\omega'$ and their
inverses.  Moreover we may assume that the base of every tassel in
$\omega$ is the base of a tassel in $\omega'$.   The inner product is
the integral with respect to uniform measure of some function of the
holonomies of the edges of $\omega'$.  These holonomies are independent
group-valued random variables distributed according to Haar measure.
Thus if any edge of $\omega'$ does not lie entirely in the range of
$\gamma$ but lies in the range of $\omega$ it will appear in the inner
product computation once, represented in some nontrivial irreducible
representation, and therefore will make the whole inner product zero.

Writing any edge $e$ of $\omega^i$ as a product of edges in $\omega'$ and
their inverses, the rightmost term in this product will be an edge of
$\omega'$ whose range lies within that of $e$ in some neighborhood of
the base of $\omega^i$.   By the previous paragraph, if the inner product of
$\Psi_\Gamma$ and $\Psi_W$ is nonzero, this edge of $\omega'$ must lie
entirely in the range of $\gamma$.  This proves the claim of the first 
paragraph, and hence the theorem.
\end{pf}

Thus $L^2(\overline{\A/\G})$ decomposes into an uncountable orthogonal
direct sum, with one countable-dimensional summand for each graph in
$M$, spanned by spin network states supported on that graph in $M$,
and one summand containing all the spin web states that are not spin
network states.

\section{The diffeomorphism-invariant Hilbert space}

The next step is to construct the diffeomorphism-invariant Hilbert
space.  Since there are very few diffeomorphism-invariant states in
$L^2(\overline{\A/\G})$, we look for diffeomorphism-invariant vectors
in a larger space.  A good choice for this larger space is the
topological dual ${\cal C}^\ast$, where ${\cal C}$ is the space of
gauge-invariant cylinder functions.  One may think of elements of
${\cal C}^\ast$ as `generalized functions' on ${\overline \A/\G}$.  We
construct diffeomorphism-invariant elements of ${\cal C}^\ast$
essentially by averaging spin network states over the action of the
diffeomorphism group, following the technique of Ashtekar,
Lewandowski, Marolf, Mour\~ao and Thiemann \cite{ALMMT}.  We also
follow their method to define an inner product on the resulting
`diffeomorphism-invariant spin network states', which allows us to
construct the diffeomorphism-invariant Hilbert space.

In general, only those diffeomorphisms of $M$ in the connected component
of the identity lift to automorphisms of the bundle $P \to M$.
However, all diffeomorphisms of $M$ lift to automorphisms of `natural'
bundles such as trivial bundles, the frame bundle, or other
bundles built from the tangent bundle using functorial constructions.
In quantum gravity it remains controversial whether one should impose
invariance under all diffeomorphisms or only those in the identity
component.  Luckily we do not need to resolve this issue here.  In what
follows, by a {\em diffeomorphism\/} we mean an element of some fixed
subgroup ${\cal D} \subseteq {\rm Diff}(M)$, all of whose elements lift
to automorphisms of $P$.  Note that with this definition all
diffeomorphisms act on $\A/\G$, $L^2(\overline{\A/\G})$, $\cal C$,
${\cal
C}^\ast$, and so on.

Given a spin network $\Gamma = (\gamma, \rho,\iota)$, the range
of $\gamma$ is a graph in $M$, say $G$.  We may write $G$ in
a unique way as a disjoint union of finitely many
points, embedded open intervals and
circles, such that none of the points has an
neighborhood in $G$ diffeomorphic to an interval embedded in $M$.  
Let ${\cal D}_\Gamma$ be the group of diffeomorphisms
mapping each of these points, intervals and circles onto itself in
an orientation-preserving way.  Let 
${\cal D}/{\cal D}_\Gamma$ be the quotient of ${\cal D}$ on the right by
${\cal D}_\Gamma$.   Note that two diffeomorphisms in the same equivalence
class of this quotient act the same way on the spin network state 
$\Psi_\Gamma$, so we can speak of the orbit $({\cal D}/{\cal D}_\Gamma) 
\Psi_\Gamma$.

\begin{proposition}\et
Given spin network states $\Psi_\Gamma,\Psi_{\Gamma'}$,
the set of elements of $({\cal D}/{\cal D}_\Gamma)\Psi_\Gamma$ having
nonzero inner product with $\Psi_{\Gamma'}$ is finite.
\end{proposition}

\begin{pf}
By Theorem 1, the inner product of $g\Psi_\Gamma$ and $\Psi_{\Gamma'}$
is zero unless $g$ takes the graph $G$ in $M$ on which $\Gamma$ is
supported to the graph $G'$ in $M$ on which $\Gamma'$ is supported.  
It follows that if we write $G$ and $G'$ as above as a union of points,
intervals and circles, $g$ establishes a one-to-one correspondence
between the points, intervals and circles of $G$ and those of $G'$.
Moreover, $[g]\in {\cal D}/{\cal D}_\Gamma$ is determined by 
this one-to-one correspondence.  Since finitely such one-to-one
correspondences are possible, there are finitely many $[g]$ 
for which $g\Psi_\Gamma$ and
$\Psi_{\Gamma'}$ have nonzero inner product.
\end{pf}

Thus it makes sense to define the quantity
\[ \bracket{\bracket{\Psi_\Gamma,\Psi_{\Gamma'}}} =
\sum_{\Phi \in ({\cal D}/{\cal D}_\Gamma)\Psi_\Gamma
}\bracket{\Phi,\Psi_{\Gamma'}} \]
for spin network states $\Psi_\Gamma$ and $\Psi_{\Gamma'}$.  We may then
extend this by sesquilinearity to all finite linear combinations of spin
network states.  It is not a priori clear that the extension is
well-defined, but in fact it is.  To see this, consider a function
$\Phi$ that can be written as a finite linear combination of spin
network states.  Consider such a decomposition, and for each graph $G$
in $M$ let $\Phi_G$ be the sum of all spin network states appearing in
the
decomposition that are supported on $G$, weighted by their coefficients.
Thus $\Psi=\sum_G \Phi_G$.  Consider another such decomposition
$\Phi=\sum_G \Phi'_G$.  By Proposition 3 the ordinary inner
product satisfies
\[ \bracket{\Phi_G,\Phi_G}=\bracket{\Phi,\Phi_G}=
\bracket{\Phi'_G,\Phi_G}=\bracket{\Phi'_G,\Phi}=
\bracket{\Phi'_G,\Phi'_G} \]
and thus $\Phi_G=\Phi'_G$.   Thus, while the exact decomposition into
spin network states is not unique, the terms $\Phi_G$ are.  But
clearly if $\Phi_G$ and $\Phi'_{G'}$ are linear combinations of
spin networks supported on the graphs $G$ and $G'$ in $M$,
respectively, then
\[  \bracket{\bracket{\Phi_G,\Phi'_{G'}}}=
\sum_i\bracket{g_i \Phi_G,\Phi'_{G'}}\]
where the $g_i$ are representatives of equivalence classes of
diffeomorphisms taking $G$ to $G'$.   This is independent of the choice
of decomposition of $\Phi_G$ and $\Phi'_{G'}$ into spin network states.
From this it follows that  $\bracket{\bracket{ \cdot,\cdot}}$ is
well-defined on finite linear combinations of spin network states.

For any spin network $\Gamma$, the linear functional
$\bracket{\bracket{ \Psi_\Gamma,\cdot }}$ extends from finite linear
combinations of spin network states to all of $\cal C$.  The space
$\cal C$ is the union over all families $C$ of the spaces of
gauge-invariant cylinder functions depending on the holonomies along
the curves in $C$.  Since each of these spaces is a Banach space
in the $\sup$ norm, one can make $\cal C$ into a topological vector
space with the inductive limit topology.  One can check that with this
topology, the {\em diffeomorphism-invariant spin network state\/}
$\bracket{\bracket{ \Psi_\Gamma,\cdot }}$ is an element of the
topological dual ${\cal C}^\ast$.  That $\bracket{\bracket{
\Psi_\Gamma,\cdot }}$ is really diffeomorphism-invariant follows from:

\begin{theorem}\et
$\bracket{\bracket{\cdot , \cdot}}$ is a positive-semidefinite,
conjugate symmetric, sesquilinear form on finite linear combinations
of spin network states.  The quotient by the null space is exactly the
quotient by the action of the diffeomorphism group.
\end{theorem}
\begin{pf}
The conjugate symmetry and sesquilinearity is obvious.  To see that it
is positive semidefinite, consider $\Phi=\sum_G \Phi_G$, with
notation as above.  We have
\[   \bracket{\bracket{\Phi,\Phi}}=
\sum_{G,G'}\bracket{\bracket{\Phi_G,\Phi_{G'}}}, \]
Note however that if $g$ is a
diffeomorphism which takes $G$ to $G',$ then $g \Phi_G$ is
supported on $G'$ and
\[ \bracket{\bracket{\Phi_G,\Phi_{G'}}} =
  \bracket{\bracket{g\Phi_G,\Phi_{G'}}}. \]
Thus, if we divide the graphs $G$ in $M$ into equivalence
classes of graphs in $M$ that are all diffeomorphic to each
other, choose a representative of each class, and choose
diffeomorphisms connecting each to the representative, we can
can replace $\Phi$ with $\Phi'=\sum_G \Phi'_G$ where now distinct
$G$ cannot be mapped to each other by diffeomorphisms, and
\[  \bracket{\bracket{\Phi,\Phi}} =\bracket{\bracket{\Phi',\Phi'}}=
\sum_G \bracket{\bracket{\Phi'_G,\Phi'_G}}.\]

To see that $\bracket{\bracket{\Phi'_G,\Phi'_G}}\geq 0$, choose
representatives $g_i$ of the equivalences classes of diffeomorphisms
which map $G$ to itself, and note that
\[  \sum_{i,j} \bracket{g_i
\Phi'_G,g_j \Phi'_G}= \sum_{i,j} \bracket{g_j^{-1}g_i
\Phi'_G, \Phi'_G}=\sum_i n_G \bracket{g_i \Phi'_G,\Phi'_G}=
n_G \bracket{\bracket{\Phi'_G,\Phi'_G}},\]
where $n_G$ is the number of diffeomorphisms $g_i$.  Since the original
inner product is positive-definite, $\bracket{\bracket{\Phi'_G,\Phi'_G}}
\geq 0$, and it is zero exactly when each $\Phi'_G$ has $\sum_i g_i
\Phi'_G = 0$.  But this condition says exactly that the $\Phi'_G$ (and
hence the original $\Phi_G$) are a sum of elements of the form $1/n_G
\sum_i (\Phi'_G -g_i \Phi'_G)$.  From this the last statement follows.
\end{pf}

If we quotient the space of finite linear combinations of spin network
states by the kernel of $\bracket{\bracket{\cdot , \cdot}}$ and then
complete it in this inner product, we obtain the {\em
diffeomorphism-invariant Hilbert space\/} $\hdiff$.  Any spin network
state $\Psi_\Gamma$ determines a diffeomorphism-invariant state
$[\Psi_\Gamma] \in \hdiff$, and also 
a continuous linear functional $\langle \langle \Psi_\Gamma,
\cdot\rangle\rangle$ on the space of cylinder functions.  The map
\[    [\Psi_\Gamma] \mapsto \langle \langle \Psi_\Gamma, \cdot\rangle\rangle
\]
extends uniquely to a continuous linear map from $\hdiff$ to ${\cal C}^\ast$.
Since this map is one-to-one, we may think of $\hdiff$ as a subspace
of the space of diffeomorphism-invariant vectors in ${\cal C}^\ast$.

\section{Problems with spin webs}

We now give examples of:
\begin{enumerate}
\item Spin web states whose inner product is nonzero but whose underlying 
webs do not have the same range.
\item A spin web state $\Psi_W$ whose orbit $({\cal D}/{\cal
D}_W)\Psi_W$ contains infinitely many distinct spin web states whose
inner product with $\Psi_W$ is nonzero.  (Here ${\cal D}_W$ is the 
set of diffeomorphisms fixing $\Psi_W$.)  
\end{enumerate}

The examples are generated out of the standard smooth function
constructed in most introductory analysis classes, whose domain and
range are $[0,1],$ which is positive on $(0,1),$ and whose value and
all order derivatives are $0$ at $0$ and $1$.  Choose one such and call
it $\alpha(x)$.  Let $\alpha_{a,b}$ be $\alpha$ composed with a linear
function so that its domain is now $[a,b],$ and let
$x_i=1/2(1+\sgn(i)(1-2^{-i})),$ an order-preserving map of the
integers into the unit interval with $0$ and $1$ as accumulation
points.  Now let
$$\beta_i^\pm= \pm 2^{-4^{|i|} } \alpha_{x_i,x_{i+1}}.$$ This
unintuitive formula describes a doubly infinite sequence of disjoint
(except for their endpoints) `blips' above and below the $x$-axis
between $0$ and $1,$ converging to both endpoints in such a fashion
that any choice of signs for each integer $i$ indicates a collection
of functions which can be pasted together to get a smooth function
on the unit interval whose graph is an embedded curve in the plane.

\begin{figure}[hbt]
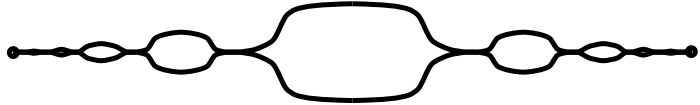

$$\pic{tassel}{40}{-20}{0}{0}$$
\label{fg:type}
\caption{The four curves $c_1,\dots,c_4$}
\end{figure}

Now let $c_1$ be the curve gotten by taking all plus signs, $c_2$ by
taking all minus signs, $c_3$ by taking a plus sign when $i$ is even and
a minus when $i$ is odd, and $c_4$ by taking a plus when $i$ is odd and
a minus when $i$ is even.  The range of these four curves is shown in 
Figure 1.

Fix a trivial $\SU(2)$ bundle over the plane.
Label each curve $c_i$ with the same representation $\rho_i$, namely
the spin-1/2 representation (i.e.\ the 2-dimensional irreducible
representation), and assign both endpoints the canonical invariant
element of $\rho_1 \tnsr \rho_2$ tensored with the canonical element
of $\rho_3 \tnsr \rho_4$, the subscripts indicating to which curve the
representation corresponds.  The family of four curves is not itself a
web, but if we use the labelings to define a function $\Psi$ of
connections in the usual way, it is easy to check by cutting each
curve in half at $x = 0$ that $\Psi$ is a spin web state.  It is also
easy to check that the holonomies of the four curves are independent
random variables with Haar measure distribution with respect to
uniform measure.

Now let $\Phi$ be defined the same way, only pick some odd $i$ and make
$c_2$ and $c_3$ take the plus rather than the minus route at $i,$ so
that $\beta_i^-$ is not in the range of the web supporting $\Phi$.  Thus
the range of the web supporting $\Phi$ is a proper subset of that for
$\Psi$.  A calculation shows that the inner product of $\Psi$ and 
$\Phi$ is nonzero.  Thus $\Phi$ gives an example of the first observation.  
In fact, the same construction gives infinitely many such $\Phi$.

For the second observation, we think of the curves $c_i$ as living in
the $xy$ plane in $\R^3$.  We consider the same $\Psi,$ and for each
$i$ consider a diffeomorphism $g_i$ which interchanges the curves
$\beta_i^+$ and $\beta_i^-$ and leaves the other $\beta_j^\pm$ fixed.
The inner product of $\Psi$ and $g_i\Psi$ is nonzero even though these
states are distinct.  Thus the spin web states $g_i\Psi$ are an
infinite class of different elements of the orbit of $\Psi$ having
nonzero inner product with $\Psi$.

Based on these examples, it would seem quite difficult to give an
effective procedure for constructing an orthonormal basis of the full
$L^2(\overline{\A/\G})$ or to give a version of `averaging over the
action of the diffeomorphism group' that would apply to spin webs that
are not spin networks.

\subsection*{Acknowledgments}

The authors would like to thank Jerzy Lewandowski and Thomas Thiemann
for helpful conversations and correspondence, and Don Marolf for
corrections to the first version of this paper.

\end{document}